\documentclass[11pt]{amsart}
\usepackage{saltscribe}
\usepackage[shortlabels]{enumitem}

\def \omij{\omega_{i, j}}
\def \TF{\mathrm{TF}}
\def \ux{\mathbf{e_1}}
\def \ut{\mathbf{e_t}}
\def \us{\mathbf{e_s}}
\def \uy{\mathbf{e_2}}

\begin{document}
\title[Universality of Polymers]{Improved universality bounds for directed polymers in the intermediate disorder regime}
\author{Pranay Agarwal}
\address{Pranay Agarwal, Department of Mathematics, University of Toronto, Toronto, Canada.}
\email{pranay.agarwal@mail.utoronto.ca}

\begin{abstract}
    We prove universality of Tracy--Widom GUE fluctuations for directed polymers in $1+1$ dimensions in the intermediate disorder regime. Building on the Lindeberg replacement method of \cite{Ran25}, we refine estimates for the measure of steep paths using probabilistic arguments. Our result extends the admissible range of the inverse temperature scaling parameter $\beta = n^{-\alpha}$ to all $\alpha \in (2/(3K+11), 1/4)$, provided the weights match $K$ moments with the log--Gamma distribution. For a general class of distributions, this gives universality for $\alpha > 2/17$, and when the third moment vanishes, this threshold improves to $\alpha > 1/10$. These results substantially broaden the known range of universality for directed polymers in the intermediate disorder regime.
\end{abstract}
\maketitle
\setcounter{tocdepth}{1}
% \tableofcontents

\section{Introduction} \label{s: intro}
The Tracy--Widom distribution is a central object in the KPZ universality class. Originally appearing in random matrix theory as the distribution of the largest eigenvalue in Gaussian ensembles, this distribution was proven to be the weak limit of the maximal energy in certain last passage problems in the seminal works of \cite{BDJ99} and \cite{Jo00}. Since then, it has been conjectured that this property should hold for a general class of last passage problems, under some moment assumptions on the random environment. Furthermore, this conjecture extends to the asymptotic free energy fluctuations of directed polymers in $1+1$ dimensions, which can be thought of as a positive-temperature version of last passage models. 

A major advancement on the conjecture for directed polymers was made in \cite{BCR13}, which rigorously established this conjecture for the exactly solvable log--Gamma polymer, though only for a restricted range of temperatures. This limitation was later overcome in \cite{KQ18}, where the result was extended to the intermediate disorder regime, in which the inverse temperature parameter $\beta$ scales with the system size $n$, i.e., $\beta = O(n^{-\alpha})$ for some $\alpha >0$. This regime first appeared in \cite{AKQ14} which proved the existence of a limiting object called the continuous directed random polymer and showed the convergence of directed polymer models to this object under appropriate scaling when $\alpha = 1/4$. Specifically, \cite{KQ18} showed that for the log--Gamma polymer, if $\beta_n \ll n^{-1/4}$, the rescaled free energy converges to the Tracy--Widom GUE distribution via an asymptotic analysis of the Fredholm determinant formula. Their proof, based on exact solvability, also extended to other models with sufficiently strong moment matching to the log--Gamma distribution by use of the perturbation method. 

Following this, \cite{Ran25} extended their findings using the Lindeberg replacement method, which relaxed the moment matching conditions and allowed for more general disorder distributions. A key part of the proof was the comparison of polymer measure restricted to a small time strip, with the Bernoulli measure on random walks. However, this comparison only holds valid when one removes paths which have a very steep slope across this strip. Let us denote them by $S$, and we shall formally define them later. In order to control the contribution from such paths they made an estimate which required $\alpha >1/5$. This restriction did not stem from the moment matching method itself but from the difficulties in estimating the measure of the set $S$. In this work we shall employ probabilistic techniques to prove that such steep paths have very low probability under the polymer measure irrespective of $\alpha$. This allows us to reduce the restriction on $\alpha$ to what is imposed by the moment matching argument. 

\subsection{Polymer models} To each vertex $(i,j) \in \bZ^2$ we assign a random weights $\omij(\beta)$ parameterized by $\beta$, the inverse temperature parameter. We define the weight of an upright path $\pi$ in $\bZ^2$ as the product of all the random weights along the path, and notate it as follows
\[
    \ell(\pi) := \prod_{v \in \pi} \omega_v(\beta).
\] 
These weights induce a probability measure on the set of paths between two fixed endpoints. The normalization constant required to define this measure is commonly referred to as the partition function.  Let $(m,n) \in \bZ_{\ge 0}^2$, we notate the partition function for paths from the origin to $(m,n)$ by
\[
    Z_{m,n}^\beta = \sum_{\pi \in \Pi_{m,n}} \ell(\pi).
\]
Here $\Pi_{m,n}$ denotes the set of all up-right paths $\pi$ which start at $(0,0)$ and end at $(m,n)$, that is the set of functions
\[
\pi: \{0,1, \dots, m+n\} \to \bZ_{\ge 0}^2
\]
such that $\pi(0)=(0,0)$, $\pi(m+n)=(m,n)$ and for each $i$, $\pi(i+1)-\pi(i)=\ux$ or $\uy$, which are the unit vectors along the $x$ and $y$ axis respectively. Similarly, we define the space and time coordinates along the vectors $\us = \uy - \ux$, and $\ut = \ux + \uy$ respectively. 

When $\beta$ is clear from context, we omit it from notation. We write the partition function as $Z_n$ when $m=n$. Occasionally we shall also deal with partition function across arbitrary endpoints $\mathbf{u \leq v} \in \bZ^2$ which we shall still denote by $Z_{\mathbf{u,v}}$ for notational simplicity.

Logarithm of the partition function is referred to as the \textit{free energy}, and the \textit{polymer measure} is formally defined as the random probability measure $\bQ_{m,n}^\beta$ on up-right paths defined by
\[
\bQ_{m, n}^\beta(\pi)=\frac{\ell(\pi)}{Z_{m,n}}.
\]
We shall use the same notational shorthands for $\bQ_{m,n}^\beta$ as we have defined for $Z_{m,n}^\beta$. For the model weights, we shall work within the following assumptions,
\begin{assumption} The model weights $(\omij)$ are independent random variables satisfying
\begin{enumerate}[1.] \label{ass}
    \item $\omega_{i,j} >0$ almost surely.
    \item $\E{\omega_{i,j}} =1$.
    \item For each $k \in \bN$, there exists constant $C_k$ such that for $\beta$ small enough,
    \[
        \bE{\abs{\omega_{i,j} -1}^k} \leq C_k \beta^k. 
    \]
    \item For each $k \in \bN$, $\E{\omij^{-k}} \to 1$ as $\beta \to 0$.
    \item There exists constants $M, C_0 > 0$, such that for all $0 < s< 1$, there exists $A(s) >0$ such that
    \[
    \P{e^{-M \beta^s} \le \omij(\beta) \le e^{M \beta^s}} \ge 1-A(s) e^{-C_0\beta^{s-1}},
    \]
    for all $\beta$ small enough.
\end{enumerate}
\end{assumption}

As discussed in \cite[Theorem 7.1]{Ran25}, these assumptions hold for a large class of interesting model weights and their mixtures. The fourth assumption does not appear in \cite{Ran25}, however all model weights that were considered in his work, satisfy this assumption. In particular, the following weight distributions come under our assumption:
\begin{enumerate}[label = \underline{Distribution \Roman{*}}, align = left]
    \item \label{ctg_1}: $\omij(\beta) = \phi_{i,j}(\beta)^{-1} e^{\beta \xi_{i,j}}$, where $\{\xi_{i,j}\}$ is a collection of independent random variables with a uniform exponential tail, and $\phi_{i,j}$ is the moment generating function of $\xi_{i,j}$.
    \item \label{ctg_2}:  $\omega_{i,j}(\beta) = (\theta -1) / X_{i,j}$, where $\theta = \theta(\beta) \sim \sigma^2/\beta^2$ as $\beta \to 0$, and $X_{i,j}$ are i.i.d.\:with $X_{i,j} \sim \mathrm{Gamma}(\theta,1)$.
\end{enumerate}

We shall improve the main results in \cite{Ran25} under these assumptions on model weights. 

\begin{theorem} \label{thm: mom_match}
    Let $\omij$ and $\omij'$ be two collections of independent random variables which satisfy Assumption \ref{ass}. Let $K\in \bN$, $K \geq 2$ be such that 
    \begin{align*}        
        \abs{\E{\omij^k} - \E{{\omij'}^k}} = \begin{cases}
            0 &\mathrm{for} \; k =1, 2,\\
            O(\beta^{K+1}) &\mathrm{for} \; k = 3, \cdots, K, 
        \end{cases}        
    \end{align*}
    for $n$ large enough. Fix $\alpha \in \br{\frac{2}{3K+11} ,\frac{1}{4}}$, then for a probability distribution $F$ on $\bR$ and a sequence of real numbers $(a_n)$ and constant $\sigma$, we have
    \begin{align}
        \frac{\log Z_n ^\beta -a_n}{ \sigma \beta^{4/3} n^{1/3}} \xrightarrow{d} F && \mathrm{if\: and\: only\: if} && \frac{\log Z_n'^\beta -a_n}{ \sigma \beta^{4/3} n^{1/3}} \xrightarrow{d} F.
    \end{align}
\end{theorem}

Weight configurations satisfying \ref{ctg_2} are referred to as the log--Gamma polymer. It is already known for this configuration that for all values of $\alpha >0$, we have Tracy--Widom weak convergence of the rescaled partition function. This in conjunction with above theorem improves the bounds stated in \cite[Theorem 2.1]{Ran25} as follows,

\begin{theorem} \label{thm: TW}
    Let $\alpha \in (2/17,1/4)$, and set $\beta=n^{-\alpha}$. Let $\xi_{i,j},\: i,j \ge 0$ be i.i.d.\,random variables with variance $\sigma^2$ and a well-defined moment generating function $\phi(t) = \E{e^{t\xi_{1,1}}}<\infty$ for all $t$ in a small neighbourhood of zero. Define
    \[
     \theta =2 + \frac{\phi(\beta)^2}{\phi(2\beta)-\phi(\beta)^2} ,\qquad\mathrm{and}\qquad a_n=2n\br{\log \phi(\beta) + \log(\theta-1) - \Psi(\theta/2)},
    \]
    where $\Psi$ is the digamma function. Let us denote the partition function between the origin and $(n,n)$ for weights $\omij = e^{-\beta \xi_{ij}}$ by $Z_n^\beta$, then we have that
    \begin{equation} \label{eq: TW}
        \frac{\log Z_n^\beta-a_n}{(4\sigma^4\beta^4 n)^{1/3}} \xrightarrow{d} \mathrm{TW}_\mathrm{GUE}.
    \end{equation}

    If $\xi_{i,j}$ satisfy the additional constraint that $\E{\xi_{i,j}^3} =0$, then \eqref{eq: TW} holds for all $\alpha \in (1/10, 1/4)$.
\end{theorem}

\subsection*{Notation} We shall frequently make use of asymptotic notations in our arguments. We say $f(n) \sim g(n)$ (resp. $f(n) << g(n)$) if $f/g \to 1$ (resp. $0$) as $n \to \infty$. Similarly, we say $f = O(g)$ if $f \leq C g$ for some arbitrary constant $C$. In general, we shall use $C, c$ to denote arbitrary positive constants which may change from equation to equation. 

We have avoided using floor and ceiling functions in many places to make notations less cumbersome. In many arguments without loss of generality we shall implicitly assume that certain fractions or fractional powers are integers. It will be easy to check that these assumptions would not affect our arguments in a non-trivial way.

\section{Fluctuations of polymer}
In this section, we establish a probabilistic upper bound on large local transversal fluctuations of paths under the polymer measure. This bound will eventually yield a control on the measure of set $S$ of steep paths. Consider a polymer model with weights following Assumption \ref{ass}, and for $0 <s <1$ define the event
\begin{equation} \label{eq: wt_bnd}
    \cW_s := \bc{e^{-M \beta^s} \le \omij (\beta) \le e^{M \beta^s} \vert i,j \in \bs{n}},
\end{equation} 
where $M$ is as defined in Assumption \ref{ass}. It follows from a simple union bound that
\[
    \P{\cW_s} \ge 1 - n^2 A(s) e^{-C_0\beta^{s-1}}.
\]
Consider an upright path $\pi$ from the origin to $(m,n)$. On the event $\cW_s$, we have that
\begin{equation} \label{eq: pth_bnd}
    e^{-2 M (n+m) \beta^s} \leq \ell(\pi) \leq e^{2 M (n+m) \beta^s}.
\end{equation}
To quantify transversal fluctuations, we measure the deviation of a path $\pi$ on the line $\cL_r := \bc{x+y = r}$ for $r \leq m+n$, as $\TF(\pi, r) := \pi(r). \us/2$. We note that after $n$ steps, the number of paths $\pi$ with $\TF(\pi, n) =k$, can be written as
\[
    a(n,k) = {n \choose n/2 +k}.
\]  
The asymptotics of $a(n,k)$ will play a key role in controlling the fluctuations of the polymer. The following estimate can be proven by direct computations. One can also refer to the proof of Theorem 1 of Chapter VII of \cite{Feller68} for details. 
\begin{lemma} \label{lemma: pnc}
    For $n$ large enough and $k << n$, we have that
    \begin{equation}
        a(n,k) =  \frac{2^{n+1/2}}{\sqrt{n \pi}} \exp \br{-(2+o(1))\frac{k^2}{n}}.
    \end{equation}
\end{lemma}
The event $\cW_s$ bounds the weight of paths, while the above lemma gives us a bound on the number of paths with large deviation. We can now combine these estimates to state a rough bound on the transversal fluctuation of the polymer:
\begin{align}    
    \bQ_n\br{\abs{\TF(\pi, n)} \geq tn^{1-s\alpha/2}} &= \frac{\sum_{\pi : \abs{\TF(\pi, n)} \geq tn^{1-s\alpha/2}} \ell(\pi)}{\sum_\pi \ell(\pi)}  \label{eq: glob_fluc}\\
    &\leq e^{4nM \beta^s} \sum_{k \geq tn^{1-s\alpha/2}} a(n,tn^{1-s\alpha/2})^2 /2^{2n} \notag\\
    &\leq \exp\br{-ct^2n^{1-s\alpha}}, \notag
\end{align}
where the last inequality follows by taking $t$ and $n$ large enough.

Our first task will be to use this upper bound on the global fluctuations of the polymer to establish an upper bound for the local fluctuations. Techniques for establishing such results are known in the literature, we shall be borrowing the proof idea from {\cite[Theorem 3]{BSS19}} (see also \cite[Proposition 2.1]{BBB25}). Similar techniques to transfer global fluctuation lower bound to a local fluctuation lower bound can be found in {\cite{Aga25}}.

\begin{proposition} \label{prop: loc_fluc}
    There exists constants $r_0, t_0$ such that for $r_0< r << n$ and $t_0 < t < r^{s\alpha/2}$, we have that 
    \[
        \bQ_n\br{\abs{\TF(\pi, r)} \geq tr^{1-s\alpha/2}} \leq \exp\br{-ct^2r^{1-s\alpha}}.
    \]
\end{proposition}
\begin{proof}
    By symmetry, it suffices to consider the event $\cA = \bc{\TF(\pi, r) \geq tr^{1-s\alpha/2}}$. Let $N = \log_2 (n/2r) -1$ and without loss of generality we assume that it is an integer. For $i=0, \cdots, N$ define the events
    \[
        \cA_i := \bc{\TF(\pi, 2^ir) \geq 2^{is\alpha/4}t(2^ir)^{1-s\alpha/2}} \cap \bc{\TF(\pi, 2^{i+1}r) \leq 2^{(i+1)s\alpha/4} t (2^{i+1}r)^{1-s\alpha/2}}.
    \]
    We now define the event $\cB = \bc{\TF(\pi, n) \geq tn^{1-s\alpha/2}}$ and claim that event $\cA$ implies either event $\cB$ or one of the event $\cA_i$ for some $i < N$. If event $\cA$ occurs, then event $\cA_i$ occur when the transversal fluctuation of the path has fallen below a certain threshold. If the fluctuations have not fallen, i.e.\,none of the events $\cA_i$ have occurred, then by our choice of $N$, we will have that the event $\cB$ must have occurred. We have depicted the argument in Figure \ref{fig:loc_fluc}. Thus,
    \begin{equation} \label{eq: prop2.2}
        \bQ_n\br{\cA} \leq \bQ_n\br{\cup_i \cA_i \cup \cB} 
                      \leq \sum_i \bQ_n\br{\cA_i} + \bQ_n\br{\cB}.  
    \end{equation}
    \begin{figure}[tb]
        \centering
        \includegraphics[width=0.6\textwidth]{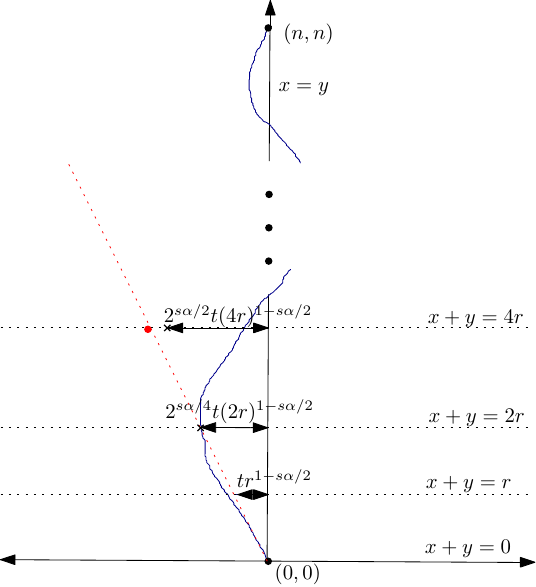}
        \caption{Depiction for the event $\cA$, showcasing event $\cA_1$ in particular. The key insight is that the fluctuations must eventually fall below the threshold; otherwise, the path exhibits a large global fluctuation.}
        % If we just analyze the journey of the polymer till the line $x+y =4r$, then on the event $\cA_1$ the polymer will be forced to have a large fluctuation on $x+y =2r$. The number of paths which can show such a large fluctuation is very small and thus the effect from the random environment becomes the dormant factor.
        \label{fig:loc_fluc}
    \end{figure}

    We have already bounded $\bQ_n(\cB)$ in \eqref{eq: glob_fluc}. To prove the required bound, it suffices to just bound the sum in the above display.

    Any path $\pi$ which falls under event $\cA_0$ will intersect the lines $\cL_r$ and $\cL_{2r}$ at points $w_1$ and $w_2$ such that 
    \[
        x_1 := w_1. \us /2r^{1-s\alpha/2} \geq t, \qquad \mathrm{and}\qquad x_2 := w_2. \us /2(2r)^{1-s\alpha/2} \leq 2^{s\alpha/4}t. 
    \]
    The idea now is to calculate the polymer measure of $\cA_0$ conditioned on $x_1$ and $x_2$ taking a particular value in this range. These calculations can be done similar to how it was done in \eqref{eq: glob_fluc}. We then maximize this conditional probability on the above range which yields an upper bound for $\cA_0$. Let us start by fixing $k_1 \geq t$, $k_2 \leq 2^{s\alpha/4}t,$ and picking points $w_1$ and $w_2$ so that $x_1 = k_1$ and $x_2 = k_2$. Define $\Pi_1$ as the set of all paths from $(0,0)$ to $w_2,$ and let $\Pi_2 \subset \Pi_1$ be the set of paths which also pass through $w_1$. We can now claim that
    \[
        \bQ_n(\cA_0 \:\vert\: x_1 = k_1, x_2 = k_2) = \frac{\sum_{\pi \in \Pi_2} \ell (\pi)}{\sum_{\pi \in \Pi_1} \ell (\pi)},
    \]
    since the contribution from the journey after $\cL_{2r}$ is the same for the numerator and the denominator. By Lemma \eqref{eq: pth_bnd}, and \ref{lemma: pnc} we can claim that the R.H.S.\ can be bounded by 
    \begin{align*}
        \frac{e^{4rM \beta^s}\abs{\Pi_2}}{\abs{\Pi_1}} &\leq e^{4rM \beta^s} \frac{a(r, k_1r^{1-s\alpha/2}) \times a(r, (2^{1-s\alpha/2}k_2-k_1)r^{1-s\alpha/2})}{a(2r, k_2(2r)^{1-s\alpha/2})} \\
        &\leq  e^{4rM \beta^s} \frac{\exp \br{ -(2 +o(1))(k_1^2+(2^{1-s\alpha/2}k_2-k_1)^2) r^{1-s\alpha}}}{\exp \br{-(2 +o(1))k_2^2 (2r)^{1-s\alpha}}} \\
        &\leq e^{4M r^{1-s\alpha}} \exp \br{ -(2+o(1))(2k_1^2 -2^{2-s\alpha/2}k_1 k_2 +2^{1-s\alpha} k_2^2)r^{1-s\alpha}}\\
        &= e^{4M r^{1-s\alpha}} \exp \br{ - (4+o(1))(k_1 - 2^{-s\alpha/2} k_2)^2r^{1-s\alpha}}.
    \end{align*}         
    Due to the restrictions on $k_1$ and $k_2$, we have that the square in the second exponent cannot be zero, in particular
    \[
        (k_1 - 2^{-s\alpha/2} k_2)^2 \geq (t - 2^{-s\alpha/4} t)^2 = c t^2.
    \] 
    Following the discussion above, we can now conclude that 
    \begin{equation*}
        \bQ_n(\cA_0) \leq  \sup_{k_1, k_2} \bQ_n(\cA_0 \:\vert\: x_1 = k_1, x_2 = k_2) \leq e^{4M r^{1-s\alpha}} e^{ct^2 r^{1-s\alpha}} \leq e^{ct^2 r^{1-s\alpha}}, 
    \end{equation*}
    for $t$ large enough. Now by simply changing the values of $t$ and $r$, we can claim that
    \begin{equation}
        \bQ_n(\cA_i) \leq \exp\br{c 2^i t^2 r^{1-s\alpha}}.
    \end{equation}
    Applying this upper bound along with \eqref{eq: glob_fluc} to \eqref{eq: prop2.2} tells us that
    \begin{equation*}
        \bQ_n\br{\cA} \leq \sum_{i=0}^N \exp\br{c 2^i t^2 r^{1-s\alpha}} + \exp\br{-ct^2n^{1-s\alpha}} \leq \exp\br{-ct^2r^{1-s\alpha}}.  \qedhere
    \end{equation*} 
\end{proof}
\subsection{Measure of paths in $S$}
Let $[a,b] \subset [0, 2n]$ be some time interval of length $n_0$. Define $S$ as the set of paths across the time interval $[a,b]$ such that 
\[
    \frac{(\pi(b) -\pi(a)). \us}{b-a} \notin \bs{-\frac{1}{2}, \frac{1}{2}} \qquad \forall \pi \in S,
\] 
and let $\tilde S$ denote the set of paths from the origin to $(n,n)$ whose restriction to $[a,b]$ is $S$. We denote the random weights of our environment by $\omij$ and define $\omij'$ as
\[
    \omij'=\begin{cases}
        \omij \qquad  i+j \notin [a,b],\\
        1\qquad \;\;\;\mathrm{otherwise}.
    \end{cases}
\]
Let us denote the polymer measure corresponding to this modified configuration as $\bQ'$. We control the total weight of steep paths under $\bQ'_n$ as follows.
\begin{theorem} \label{thm: S}
    Let $\alpha \in (0,1/4), s \in (0,1)$. Then for all $\gamma >0$, $n_0, a, 2n -b \geq n^\gamma$, and $ n_0 << n$, we have that on the event $\cW_s$,
    \begin{equation} \label{eq: thm_S}
        \bQ'_n \br{\tilde S} \leq \exp \br{-cn^{\gamma(1-s\alpha)}}. 
    \end{equation}
\end{theorem}
\begin{proof}
    By symmetry of the model about the line $x+y =n$, we can assume that $a \leq n$ without any loss of generality. Let us define the event
    \[
        \cA = \bc{\pi \in \Pi_n \: :\: \abs{\TF(\pi, a)} \leq 2t_0 a^{1-s\alpha/2}},
    \]
    where $t_0$ is as defined in Proposition \ref{prop: loc_fluc}. If $\omij$ satisfy Assumption \ref{ass}, then so does $\omij'$. Therefore, by Proposition \ref{prop: loc_fluc} we have that
    \begin{equation} \label{eq: S1} 
        \bQ_n'(\tilde S) \leq \bQ_n'(\tilde S \cap \cA) + \bQ_n'(\cA^c) \leq \bQ_n'(\tilde S \cap \cA) + \exp \br{-cn^{\gamma(1-s\alpha)}}.
    \end{equation}
    Thus, it suffices to show that $\bQ_n'(\tilde S \cap \cA)$ satisfies the bound in \eqref{eq: thm_S}. Let 
    \begin{align*}
        I_a = \bc{u \in \cL_a: u.\us \in (-2ta^{1-s\alpha/2}, 2ta^{1-s\alpha/2})},
    \end{align*}
    denote the set of valid starting points on $\cL_a$ on the event $\cA$. Similarly, for a fixed starting point $u \in I_a$, let  
    $$I_b^u =  \bc{v \in \cL_b: (v-u).\us \notin (-n_0/2, n_0/2)}$$
    denote the set of endpoints on $\cL_b$ which can be seen on the event $S$. Thus, we can say that
    \begin{align}
        \bQ_n'(\tilde S \cap \cA) &= \sum_{u \in I_a} \bQ_n'(\pi(a) =u,  \pi(b) \in I_b^u) \notag \\
        &\leq \sum_{u \in I_a} \frac{\sum_{v \in I_b^u} Z_{(0,0), u}. Z_{u,v} . Z_{v, (n,n)}}{Z_{(0,0), u}.Z_{u, (n,n) }} \notag \\
        &\leq \sum_{u \in I_a} \frac{\sum_{v \in I_b^u} Z_{u,v} . Z_{v, (n,n)}}{Z_{u, (n,n) }}  = \sum_{u \in I_a} \bQ_{u, (n,n)}'(\abs{\TF(\pi, b)} \geq n_0/2)  \label{eq: S1}.
        % \stackrel{d}{=} \sum_{u \in I_a} \bQ_{(0,0), (n,n)- u}'(\abs{\TF(\pi, n_0)} \geq n_0/2).
    \end{align}
    Now by the same arguments as Proposition \ref{prop: loc_fluc}, we can claim that for each $u \in I_a$,
    \[
        \bQ_{u, (n,n)}'(\abs{\TF(\pi, b)} \geq n_0/2) \leq \exp \br{-cn_0^{(1-s\alpha)}} \leq \exp \br{-cn^{\gamma(1-s\alpha)}}.
    \]
    Applying the above bound to \eqref{eq: S1}, we get that
    \begin{equation}
        \bQ_n'(\tilde S \cap \cA) \leq \abs{I_a}\exp \br{-cn^{\gamma(1-s\alpha)}} \leq  \exp \br{-cn^{\gamma(1-s\alpha)}}.   
    \end{equation}
    The last inequality is a consequence of $\abs{I_a} << n$. The proof now follows from \eqref{eq: S1}.
\end{proof}

\section{Proof of Main Theorems}
We are now ready to present the proof of Theorem \ref{thm: mom_match}. The bulk of the argument will be the same as in \cite[Section 6]{Ran25}. However, we shall add some extra arguments to control the error from points very close to the endpoints. This will allow us to use Theorem \ref{thm: S} and remove the $\alpha >1/5$ restriction.

\begin{proof}[Proof of Theorem \ref{thm: mom_match}]
    Let us assume that
    \begin{align*}
        \frac{\log Z_n ^\beta -a_n}{ \sigma \beta^{4/3} n^{1/3}} \xrightarrow{d} F,
    \end{align*}
    where all symbols are as defined in Theorem \ref{thm: mom_match}. This is equivalent to saying that for any smooth function $f$ whose first $K+1$ derivatives are bounded, we have that 
    \[
        \lim_{n \to \infty} \E{f\br{\frac{\log Z_n ^\beta -a_n}{ \sigma \beta^{4/3} n^{1/3}}}} = \E{f(X)},
    \]
    for $X$ following the distribution $F$. Thus, it suffices to prove that 
    \begin{equation} \label{eq: mom_mat}
        \lim_{n \to \infty} \E{f\br{\frac{\log Z_n ^\beta -a_n}{ \sigma \beta^{4/3} n^{1/3}}}} - \E{f\br{\frac{\log Z_n'^\beta -a_n}{ \sigma \beta^{4/3} n^{1/3}}}} = 0.
    \end{equation}
    We replace each weight $\omij$ with its counterpart $\omij',$ one at a time, to get from the environment defined $\omij$ to the one defined by $\omij'$. As such, at any step we will have a mix of $\omega$ and $\omega'$ in our environment. To keep notations simple, we shall always denote the current configuration of random weights by $\varpi_{i,j}$. Now let $v \in \bZ^2_{\geq 0}$ be a vertex whose weight has not been replaced, and define
    \begin{align*}
        W_n(v) = \sum_{v \in \pi} \prod_{i : \pi(i) \neq v} \varpi_{\pi(i)}, &&\mathrm{and}&&
        V_n(v) = \sum_{\pi \in \Pi_n : v \notin \pi} \prod_{i=1}^{2n} \varpi_{\pi(i)}.
    \end{align*}
    Following the Linderberg argument in \cite[Page 6]{Ran25}, we have that
    \begin{align} \label{eq: linderberg}
        \zeta_v &:=\bE \abs{f\br{\frac{\log(V_n + \omega_v W_n) -a_n}{ \sigma \beta^{4/3} n^{1/3}}} - f\br{\frac{\log(V_n + \omega_v' W_n) -a_n}{ \sigma \beta^{4/3} n^{1/3}}}} \notag\\
        &\leq \frac{C \beta^{K+1}}{\sigma \beta^{4/3} n^{1/3}} \br{\sum_{k=3}^{K+1} \E{\frac{W_n^k}{(V_n +W_n)^k} +\frac{W_n^{K+1}}{(V_n + \omega_v W_n)^{K+1}}+\frac{W_n^{K+1}}{(V_n + \omega_v' W_n)^{K+1}}}}. 
    \end{align}
    The above inequality is a bound on the error for replacing one random weight in our environment. Using Assumption \ref{ass}.4 and positivity of $V_n, W_n$ and $\omega_v$, we can now claim that 
    \[
        \E{\frac{W_n^{K+1}}{(V_n + \omega_v W_n)^{K+1}}} = \E{\omega_v^{-(K+1)} \frac{(\omega_v W_n)^{K+1}}{(V_n + \omega_v W_n)^{K+1}}} \leq \E{\omega_v^{-(K+1)}} \leq 2,
    \]
    for $\beta$ small enough. Thus, we can naively estimate that the summands in the right-hand side of \eqref{eq: linderberg} are each less than 2, which gives us that
    \begin{equation} \label{eq: endpoints}
        \zeta_v \leq \frac{C (K+1) \beta^{K+1}}{\sigma \beta^{4/3} n^{1/3}} \implies \sum_{v \in \Delta_0} \zeta_v \to 0,
    \end{equation}
    where $\Delta_0 = \bc{u \in \bZ^2_{\geq 0} : u. \ut \leq n^{1/6}, \;\mathrm{or}\; u.\ut \geq 2n - n^{1/6}}$, denotes a set of points which are very close to the endpoints. Thus, the above approximation allows us to ignore the contribution from vertices close to the endpoint.

    Now let $\delta >0$ be small, and define $n_0 = \beta^{-4(1-\delta)}$. With this choice of $n_0$, we have that $\beta << n_0^{-1/4}$ and thus the polymer measure will behave similar to the random walk measure on strips of width $n_0$. The exact details and the intuition behind this idea are discussed in \cite[Section 6]{Ran25}.

    We now divide the time interval $(n^{1/6}, 2n -n^{1/6})$ into subintervals of length $n_0$, and denote the set of all points falling in the $i^\mathrm{th}$ subinterval by $\Delta_i$, for $i \geq 1$. Since each of the strips $\Delta_i$ are away from our endpoints, we can apply Theorem \ref{thm: S} to these strips. This, in conjunction with the Bernoulli bridge comparison in \cite[Proof of Theorem 2.2]{Ran25} gives us that for each $i\geq 1$ and $\alpha < 1/4$,
    \[
    \sum_{v \in \Delta_i} \zeta_v \leq \frac{C (K+1) \beta^{K+1} \log n}{\sigma \beta^{4/3} n^{1/3}}.
    \]
    Using the above inequality along with \eqref{eq: endpoints} allows us to bound the total error as follows,
    \[
        \E{f\br{\frac{\log Z_n ^\beta -a_n}{ \sigma \beta^{4/3} n^{1/3}}}} - \E{f\br{\frac{\log Z_n'^\beta -a_n}{ \sigma \beta^{4/3} n^{1/3}}}} = \frac{C (K+1) \beta^{K+1} \log n}{\sigma \beta^{4/3} n^{1/3}} . \frac{n}{n_0} = O(n^\lambda \log n),
    \]  
    where
    \[
        \lambda = \br{(K+1) - \frac{4}{3} + 4(1-\delta)} \alpha - \frac{2}{3}.
    \]
    Now, for any $\alpha > 2/(3K+11)$ we can choose a $\delta > 0$ so that $\lambda <0$. This proves \eqref{eq: mom_mat} and the desired result follows.
\end{proof}

We now recall the statement for the Tracy--Widom convergence of the partition function in the log--Gamma polymer.

\begin{theorem}[{\cite[Theorem 2.1]{KQ18}}] \label{thm: TW-log_G}
Let $\theta = 1/\beta^2$ and consider the log--Gamma polymer with weights $\omega_{i,j}^{-1} \sim \mathrm{Gamma}(\theta,1)$, and set $\beta = n^{-\alpha}$ for some $\alpha < 1/4$. Then, as $n \to \infty$,
\[
\frac{\log Z_n + 2n \Psi(\theta/2)}{-(\Psi''(\theta/2) n)^{1/3}} \xrightarrow{d} \mathrm{TW}_{\mathrm{GUE}},
\]
where $\Psi$ is the digamma function, defined as the derivative of the log--Gamma function, and $\mathrm{TW}_{\mathrm{GUE}}$ denotes the Tracy--Widom GUE distribution.
\end{theorem}

Let us now present the proof of Theorem \ref{thm: TW}. The overall argument is the same as the one presented in \cite[Section 7]{Ran25} in principle. However, there are some extra checks that we will be doing along the proof. In particular, we shall show that Assumption \ref{ass}.4 is also satisfied by \ref{ctg_1} and \ref{ctg_2}, apply our version of Theorem \ref{thm: mom_match}, and provide a sufficient condition for Tracy--Widom convergence for $\alpha$ as small as $1/10$. 

\begin{proof}[Proof of Theorem \ref{thm: TW}]
    Let $\xi_{ij}$ be as defined in the statement of the theorem and define $\omij = e^{\beta \xi_{ij}} \phi(\beta)^{-1}$. We now define $\omij'$ as a collection of weights satisfying \ref{ctg_2}, i.e.\ $\omij' = (\theta-1) /X_{ij}$ with $X_{ij} ~ \mathrm{Gamma}(\theta,1)$, where 
    \begin{equation} \label{eq: theta_est}
        \theta =2 + \frac{\phi(\beta)^2}{\phi(2\beta)-\phi(\beta)^2} \sim \frac{1}{\sigma^2 \beta^2}.
    \end{equation}
    From the definitions of $\omij$ and $\omij'$, their moments are given by 
    \begin{align} \label{eq: mom}
        \E{\omij^k} = \phi(k \beta) \phi(\beta)^{-k}, &&\mathrm{and} &&\E{\omij'^k} = (\theta -1)^k \frac{\Gamma(\theta -k)}{\Gamma(\theta)}. 
    \end{align}
    From here, it is easy to verify that $\omij$ and $\omij'$ satisfy Assumption \ref{ass}.4, and the rest of the assumptions have already been proven to be true in \cite[Section 7]{Ran25}.

    With the choice of $\theta$ in \eqref{eq: glob_fluc}, we have that $\E{\omij^k} = \E{\omij'^k}$ for $k =1,2$. This is a routine check, see \cite[Theorem 2.1]{Ran25} for details. Thus, the two collections of random variables satisfy the hypothesis of Theorem \ref{thm: mom_match} with $K =2$. Furthermore, we can apply Theorem \ref{thm: TW-log_G} to the partition function $Z_n'$ defined by the weights $\omij'$ to claim that
    \begin{equation*}
        \frac{\log Z_n' - 2n (\log (\theta-1)-\Psi(\theta/2))}{-(\Psi''(\theta/2) n)^{1/3}} \xrightarrow{d} \mathrm{TW}_{\mathrm{GUE}}.
    \end{equation*}
    Applying Theorem \ref{thm: mom_match}, tells us that for $\alpha \in(2/17, 1/4)$, we have that the partition function $Z_n^\beta$ (as defined in the statement of the theorem) satisfies
    \begin{equation} \label{eq: TW_conv}
        \frac{\log Z_n^\beta - 2n (\log \phi(\beta) + \log (\theta-1) -\Psi(\theta/2))}{-(\Psi''(\theta/2) n)^{1/3}} \xrightarrow{d} \mathrm{TW}_{\mathrm{GUE}}.
    \end{equation}
    Performing a Taylor expansion tells us that the digamma function can be expressed as
    \begin{align*}    
        \Psi(\theta) = \log \theta - \frac{1}{2\theta} - \frac{1}{12\theta^2} + O\left(\frac{1}{\theta^4}\right), \quad\mathrm{and} &&
        \Psi''(\theta) = -\frac{1}{\theta^2} - \frac{1}{\theta^3} - \frac{1}{2\theta^4} + O\left(\frac{1}{\theta^5}\right).
    \end{align*}
    Using the above approximation along with \eqref{eq: theta_est}, allows us to simplify the denominator of \eqref{eq: TW_conv} as follows
    \[
        \frac{\log Z_n^\beta- 2n (\log \phi(\beta) + \log (\theta-1) -\Psi(\theta/2))}{(4\sigma^4\beta^4 n)^{1/3}} \xrightarrow{d} \mathrm{TW}_\mathrm{GUE}.
    \]
    This proves the first part of our claim. We would now like to show that we can increase the range of $\alpha$ to $(1/10, 1/4)$ under the additional assumption that $\E{\xi_{i,j}^3} =0$. In light of the above argument, it suffices to show that for the same choice of $\theta$, $\omij$ and $\omij'$ satisfy the assumption of Theorem \ref{thm: mom_match} for $K=3$. Thus, we show that for $n$ large enough,
    \begin{equation} \label{eq: 3rd}
        \abs{\E{\omij^3} - \E{\omij'^3}} \leq C \beta^4.
    \end{equation}    
    The moments of $\omij/\omij'$ can be expressed in terms of the moment generating function $\phi$, thus we start by expressing $\phi$ in terms of the moments of $\xi_{ij}$. We also note how the expression changes when $\phi$ is raised to a power $k \in \bZ$, 
    \begin{align} \label{eq: phi_exp}
        \phi(k\beta) = 1 + \sigma^2 (k\beta)^2 /2 + O(\beta^4), &&\mathrm{and}&& \phi(\beta)^{k} = 1 + k \sigma^2 \beta^2 /2 + O(\beta^4).
    \end{align}
    Combining the above estimates with \eqref{eq: mom}, gives us that
    \begin{equation} \label{eq: omij}
        \E{\omij^3} = 1 + 3 \sigma^2 \beta^2 + O(\beta^4). 
    \end{equation} 
    Similarly, we shall now approximate the third moment of $\omij'$,
    \begin{equation*}
        \E{\omij'^3} =  \frac{(\theta-1)^2}{(\theta-2)(\theta-3)} = \frac{\phi(2 \beta)^2}{\phi(\beta)^2 (2\phi(\beta)^2 - \phi(2 \beta))}.
    \end{equation*}
    A quick calculation using \eqref{eq: phi_exp} shows that 
    \[
      2\phi(\beta)^2 - \phi(2 \beta) = 1+ O(\beta^4),  
    \]
    and the overall expression simplifies to 
    \[
        \E{\omij'^3} = 1 + 3 \sigma^2 \beta^2 + O(\beta^4).
    \]
    Comparing the above estimate with \eqref{eq: omij} tells us that \eqref{eq: 3rd} holds and completes the proof.
\end{proof}

\subsection*{Acknowledgement} 
I am grateful to B\'alint Vir\'ag for invaluable guidance and insightful discussions throughout this project.

\bibliography{bib}
\bibliographystyle{alpha}
\end{document}